\documentclass[twoside,11pt]{article}

\usepackage{jmlr2e}
\usepackage{amsmath}

\usepackage{mathtools}
\usepackage[normalem]{ulem} 
\usepackage{subcaption}
\usepackage{enumerate}
\usepackage[ruled,vlined]{algorithm2e}
\usepackage{amsfonts}
\usepackage{amssymb}
\usepackage{booktabs}
\usepackage{multirow}
\usepackage{mathrsfs}
\usepackage{hyperref}
\usepackage{listings}

\newcommand{\remove}[1]{{}}
\newcommand{\cut}[1]{}

\RequirePackage[normalem]{ulem} 
\RequirePackage{color}\definecolor{RED}{rgb}{1,0,0}\definecolor{BLUE}{rgb}{0,0,1} 

\newcommand{\RR}{\mathbb{R}}

\newcommand{\prox}{\mathbf{prox}}

\DeclareMathOperator*{\Min}{minimize}

\newcommand{\bc}{\begin{center}}
\newcommand{\ec}{\end{center}}

\newcommand{\bdm}{\begin{displaymath}}
\newcommand{\edm}{\end{displaymath}}

\newcommand{\beq}{\begin{equation}}
\newcommand{\eeq}{\end{equation}}

\newcommand{\bfl}{\begin{flushleft}}
\newcommand{\efl}{\end{flushleft}}

\newcommand{\bt}{\begin{tabbing}}
\newcommand{\et}{\end{tabbing}}

\newcommand{\beqn}{\begin{align}}
\newcommand{\eeqn}{\end{align}}

\newcommand{\beqs}{\begin{align*}} 
\newcommand{\eeqs}{\end{align*}}  

\usepackage{color}

\definecolor{myred}{rgb}{0.6,0,0}
\definecolor{mygreen}{rgb}{0,0.6,0}
\definecolor{mygray}{rgb}{0.5,0.5,0.5}
\definecolor{mymauve}{rgb}{0.58,0,0.82}

\SetCommentSty{mycommfont}

\newcommand{\pkg}{{TMAC}}
\newcommand{\cpp}{{C++11}}
\newcommand{\repo}{\url{https://github.com/uclaopt/TMAC}}

\jmlrheading{1}{2016}{xx-xx}{x/xx}{xx/xx}{Brent Edmunds, Zhimin Peng and Wotao Yin}

\ShortHeadings{TMAC}{Brent Edmunds, Zhimin Peng and Wotao Yin}
\firstpageno{1}

\begin{document}

\title{\pkg: A Toolbox of Modern Async-Parallel, Coordinate, Splitting, and Stochastic Methods}

\author{\name Brent Edmunds \email brent.edmunds@math.ucla.edu
       \AND
       \name Zhimin Peng \email zhimin.peng@math.ucla.edu
       \AND
        \name Wotao Yin \email wotaoyin@math.ucla.edu 
       \AND
       \addr Department of Mathematics\\
       University of California, Los Angeles\\
       Los Angeles, CA 90095, USA}

\maketitle

\begin{abstract}
\pkg~is a toolbox written in C++11 that implements algorithms based on a set of modern methods for large-scale optimization. It covers a variety of optimization problems, which can be both smooth and nonsmooth, convex and nonconvex, as well as constrained and unconstrained. The algorithms implemented in \pkg, such as the coordinate update method and operator splitting method, are scalable as they decompose a problem into simple subproblems. These algorithms can run in a multi-threaded fashion, either synchronously or asynchronously, to take advantages of all the cores available. \pkg~architecture mimics how a scientist writes down an optimization algorithm. Therefore, it is easy for one to obtain a new algorithm by making simple modifications such as adding a new operator and adding a new splitting, while maintaining the multicore parallelism and other features. The package is available at \repo.
\end{abstract}


\section{Introduction}
\pkg~is a toolbox for optimization that implements algorithms based on a set of modern methods for large-scale optimization. The toolbox covers a variety of optimization problems, which can  be both smooth and nonsmooth, convex and nonconvex, as well as constrained and unconstrained.  \pkg~is designed for fast prototyping of scalable algorithms, which can be  single-threaded or multi-threaded, and the multi-threaded code can run either synchronously or asynchronously.

Specifically, \pkg~implements  algorithms based on  the following methods:
\begin{itemize}
\item \textbf{Operator splitting}: a collection of methods that decompose problems into simple subproblems. The original problem often takes the following forms: (i) minimizing $f_1(x)+\cdots+f_n(x)$, (ii) finding a solution $x$ to $0\in T_1(x)+\cdots +T_n(x)$, and (iii) minimizing $f_1(x_1)+\cdots+f_n(x_n)$ subject to linear constraints $A_1 x_1+\cdots A_m x_m=b$. In addition, any function $f_i$ can compose with a linear operator, e.g. $f(x) = g(Ax)$.

\item \textbf{Coordinate update}: a collection of methods that find a solution $x$ by updating one, or a few, of its elements each time. The coordinate ordering can follow the random, cyclic, shuffled cyclic, and greedy rules.

\item \textbf{Parallel} 
coordinate updates (either synchronous or asynchronous).
\end{itemize}
These methods are reviewed in Section~\ref{sec:review} below.

\pkg~is not a modeling language, but an algorithm development toolbox. 
The usage of this toolbox is demonstrated in the following examples:
 linear system of equations;
quadratic programming;
empirical risk minimization (e.g. $\ell_1$ and $\ell_2$ regularized  regression);
support vector machine;
portfolio optimization;
and nonnegative matrix factorization. \pkg~can use multiple cores efficiently to solve these  problems because it exploits their underlying structures. 

\subsection{Coding and design}
\pkg~leverages the \cpp~standard\footnote{C++ is standardized by ISO (The International Standards Organization)  The original C++ standard was issued in 1998. A major update to the standard, C++11, was issued in 2011.} and object-oriented design, striking for efficiency, portability, and code readability. The package is written in C++ as Matlab does not currently support shared memory programming. The thread library, a new feature of the C++11 standard, provides multithreading that is invariant to operating system. \pkg~can be compiled by C++ compilers under Linux, Mac, and Windows. We design the package so that the user can implement a sophisticated operator-splitting, coordinate-update, and sequential or parallel algorithm  with little effort. Our codebase is separated into layers: executables, multicore drivers, schemes, operators, and numerical linear algebra\footnote{For the best performance, BLAS is called for numerical algebra operations (e.g., the product of a matrix and a vector plus another vector). While \pkg~can parallelize coordinate updates, it is also possible to parallelize numerical algebra operations by linking \pkg~with a parallel BLAS package such as ScaLAPACK\citep{blackford1997scalapack}.}  that correspond to different algorithmic components. \pkg~is used by combining objects from each layer.
As a result, our codes  are short, clean, and thus easy to read and modify.




\subsection{Download and installation}
The \pkg~package can be accessed from GitHub at~\repo. The package runs on Linux, Mac, and Windows operating systems.
\subsection{Literature}\label{sec:review}
\subsection*{Operator splitting methods:}  These methods solve complicated optimization and monotone inclusion problems by simple subproblems. They started to appear in the 1950s for solving partial differential equations and feasibility problems and were rapidly developed during the 1960s--1980s. Several splitting methods such as Forward-Backward \citep{Passty1979_ergodic}, Douglas-Rachford \citep{DR56} (which is equivalent to ADMM \citep{GabayMercier1976_dual,GlowinskiMarroco1975_approximation}), and Peaceman-Rachford \citep{PR} were introduced. Recently,
operator splitting methods such as ADMM and Split Bregman \citep{GoldsteinOsher2009_split} (also see \citep{WangYangYinZhang2008_NewAlternating}) have found new applications in image processing,
statistical and machine learning, compressive sensing, and control. New methods such as primal-dual splitting \citep{Condat2013_primaldual,Vu2013_splitting}, three-operator splitting \citep{DavisYin2015_threeoperator}, and other primal-dual splitting methods \citep{LiShenXuZhang2015_multistep,ChenHuangZhang2016_primaldual,ChenHuangZhang2016_primalduala} have appeared, and they are designed to solve more complicated problems.


\subsection*{Coordinate update methods:}
As the name suggests, these methods update the selected one, or a few, elements of the variable at each iteration. The original  coordinate descent method \citep{Hildreth1957_quadratic,Warga1963_minimizing,SargentSebastian1973_convergence,LuoTseng1992_convergence} developed in 1950s and analyzed in the 1960s--1990s minimizes the original objective function with respect to the selected coordinates. Later developments such as \citep{GrippoSciandrone2000_convergence,TsengYun2009_coordinate,TsengYun2009_blockcoordinate,XuYin2013_block,BolteSabachTeboulle2014_proximal} have allowed surrogates for the objective function that are often easier or more efficient to minimize. Lately, coordinate descent has been extended so that each update no longer minimizes a function, but instead the update applies an operator, such as the coordinate projection of an operator or a coordinate-wise fixed-point to an operator \citep{CombettesPesquet2015_stochastic,BianchiHachemFranck2014_stochastic,PengXuYanYin2015_arock,PengWuXuYanYin2016_coordinate}. Hence, we call it coordinate update instead of coordinate descent.

The initial coordinate selection rule  is cyclic selection. It was widely used before other rules such as random \citep{Nesterov2012_efficiency,RichtarikTakac2014_iteration,LuXiao2015_complexity}, shuffled cyclic, greedy \citep{BertsekasBertsekas1999_nonlinear,LiOsher2009_coordinate,TsengYun2009_coordinate,PengYanYin2013_parallel,NutiniSchmidtLaradjiFriedlanderKoepke2015_coordinate}, and parallel \citep{BradleyKyrolaBicksonGuestrin2011_parallel,richtarik2016parallel} started to appear and gain popularity.

\subsection*{Asynchronous parallel methods:}

In  parallel algorithms, multiple agents attempt to solve a problem. The agents, to solve the problem, must exchange data. An algorithm is synchronous if all the agents must finish computing before they exchange data, and only after the exchange is completed can they start the next computing cycle.  Synchronization requires every agent to wait for the slowest agent
(or the one solving the most difficult subproblem) to finish computing before communicating. On shared memory architectures\footnote{On multicore architectures, agents can be threads or processes. Threads automatically share memory, whereas processes require inter-process communication.}, the synchronization of communication leads to bus contention. The agents in an asynchronous parallel algorithm, however, can run continuously; they can compute with whatever information they have, even if the latest information from other agents
has not arrived; they write their results to the shared memory while other agents are still computing. Async-parallel methods can be traced back to~\citep{chazan1969chaotic} for solving systems of linear equations.
For function minimization,~\citep{bertsekas1989parallel} introduced an async-parallel gradient projection method. Convergence rates are obtained in~\citep{tseng1991rate-asyn}.

For fixed-point problems, async-parallel methods date back to~\citep{baudet1978asynchronous}. In  the pre-2010 methods \citep{bertsekas1983distributed,BMR1997asyn-multisplit,el1998flexible,Baz200591} and the review~\citep{frommer2000asynchronous}, each agent updates its own subset of coordinates. Convergence is established under the \emph{$P$-contraction} condition and its variants~\citep{bertsekas1983distributed}. Recently, the works \citep{nedic2001distributed,recht2011hogwild,liu2013asynchronous,liu2014asynchronous,hsieh2015passcode} introduced async-parallel stochastic methods for function minimization.
For fixed-point problems,~\citep{PengXuYanYin2015_arock} introduced async-parallel stochastic methods (ARock), as well as several applications in optimization.


\subsection*{Software packages:}
There exist several packages for solving optimization problems based on splitting and coordinate methods. 
TOFCS~\citep{becker2011templates} is a framework for solving convex cone problems with first-order methods. 
SCS \citep{ocpb:16} is a convex cone program solver that applies operator splitting methods to an equivalent feasibility problem. 
Epsilon \citep{wytock2015convex} is a package for solving general convex programming by using fast linear and proximal operators.
Though these packages solve convex optimization problems from various applications, they are sequential and do not take advantage of
the multicore systems.  
PASSCoDe \citep{hsieh2015passcode} implements the muticore parallel dual coordinate descent method that solves $\ell_2$ regularized 
empirical risk minimization problems.
APPROX\citep{fercoq2015accelerated} implements  parallel coordinate descent, stochastic dual ascent, and accelerated gradient descent for block seperable objective functions.
AC-DC\citep{richtarik2016parallel} contains a suite of serial, parallel, and distributed coordinate descent algorithms for LASSO, Elastic Net, SVM, and sparse SVM.

\section{Case study}\label{sec:quick_start}
To illustrate the  usage of \pkg, consider the $\ell_1$ regularized logistic regression problem~\citep{ng2004feature} .
\begin{equation}\label{eq:l1_log}
\Min_{x \in \RR^n} \lambda \|x\|_1 + \sum_{i = 1}^m \log\left(1 + e^{-b_i \cdot a_i^T x}\right),
\end{equation}
where $\{(a_i, b_i)\}_{i = 1}^m, (a_i \in \RR^n, b \in \{1, -1\})$ is the training dataset. For demonstration purposes, we simply set the regularization parameter
$\lambda$ to 1, and the maximum number of epochs to 100, and use \pkg~to solve~\eqref{eq:l1_log} on a machine with 64GB of memory and two Intel Xeon E5-2690 v2 processors (20 cores in total). The following are the \texttt{tmac\_fbs\_l1\_log} commands to solve the model
on the news20 dataset\footnote{ This dataset is from \url{http://www.csie.ntu.edu.tw/~cjlin/libsvmtools/datasets}. In this case, $m=19,996$ and $n=1,355,191$. } with 1 thread, 4 threads, and 16 threads.
\begin{lstlisting}[language=bash]
# ------------------- running with 1 thread -----------------------#
$ tmac_fbs_l1_log -data news20.svm -epoch 100 -lambda 1 -nthread 1
[some output skipped]
Computing time  is: 29.53(s).
# ------------------- running with 4 threads ----------------------#
$ tmac_fbs_l1_log -data news20.svm -epoch 100 -lambda 1 -nthread 4
[some output skipped]
Computing time  is: 11.01(s).
# ------------------- running with 16 threads ---------------------#
$ tmac_fbs_l1_log -data news20.svm -epoch 100 -lambda 1 -nthread 16
[some output skipped]
Computing time  is: 3.87(s).
\end{lstlisting}
The flags \texttt{-data, -epoch, -nthread, -lambda} are  for the data file, maximum number of epochs,
total number of threads, and regularization parameter $\lambda$ respectively. We can see that the command-line
tool is easy to use. Beyond the simplicity, \pkg~is also efficient in the sense that the solving time is
less than 30 seconds for a problem with more than 1 million variables. We can observe that using 16 threads can achieve approximately 8 times of speedup, reducing the run time to under 4 seconds. Next, we show the major components of the source
codes for building \texttt{tmac\_fbs\_l1\_log}.

We solve \eqref{eq:l1_log} with the forward-backward splitting scheme
\begin{equation}\label{eq:fbs_l1_log}
  x^{k+1} = \underbrace{\underbrace{\prox_{ \eta \lambda \|\cdot\|_1}}_{\text{backward operator}}
  \bigg(\underbrace{x^k - \eta \, \nabla_{x} \big(\sum_{i = 1}^m \log (1 + e^{-b_i \cdot a_i^T x^k}}_{\text{forward operator}})\big)\bigg)}_{\text{forward-backward splitting scheme}},
\end{equation}
where the gradient step of logistic loss and the proximal of $\ell_1$ norm correspond to the
forward operator\footnote{A forward operator computes a (sub)gradient at the current point and takes a negative (sub)gradient step to obtain a new point.} and backward operator\footnote{A backward operator typically solves an optimization problem, and its optimality condition yields a (sub)gradient taken at the new point.} respectively. Algorithm~\ref{alg:fbs_l1_log} shows the details of
implementing~\eqref{eq:fbs_l1_log} in an asynchronous parallel coordinate update fashion.

\begin{algorithm}[H]\label{alg:fbs_l1_log}
\DontPrintSemicolon
  \SetKwInOut{Input}{Input}\SetKwInOut{Output}{output}
  \Input{$A, b$ and $x$ are shared variables, $p$ agents, $K > 0$.}
  \textbf{Initialization:} \\
  \quad $\text{foward}(x) := x - \eta \, \nabla_x \,\sum_{i = 1}^m \log (1 + e^{-b_i \cdot a_i^T x})$  \tcp*{forward operator}
  \quad $\text{backward}(x) := \prox_{\eta \lambda  \|\cdot\|_1} (x)$ \tcp*{backward operator}
  \quad $\text{fbs}(x) := \text{backward}(\text{forward}(x))$ \tcp*{foward-backward splitting scheme}
  \quad create $p$ computing agents \\
  \While{each of the $p$ agents continuously}{
    \textbf{selects} $i \in \{1, ..., n\}$ based on some index rule \;
    \textbf{updates} $x_i \gets x_i - \eta \left(x_i - \text{fbs}_i (x)\right)$
  }
  \caption{\pkg~for $\ell_1$ logistic regression.}
\end{algorithm}
The snippet of code (extracted from
\texttt{apps/tmac\_fbs\_l1\_log.cc}) shown in Listing \ref{code:l1log} implements Algorithm~\ref{alg:fbs_l1_log} with the \pkg~package.
Specifically, line 3 defines \texttt{forward} as an operator of type \texttt{forward\_grad\_for\_log\_loss<SpMat>}
initialized by the pointers to the data \texttt{A}, which is a sparse matrix, and label \texttt{b}, which is a dense vector. Line 5  defines a
\texttt{prox\_l1} operator (\texttt{backward}) initialized by $\lambda$ and step size $\eta$.
Line 7 defines a forward-backward splitting scheme (\texttt{fbs}) with the previously defined forward operator,
backward operator and the address of the unknown variable \texttt{x}. Line 9 calls the multicore driver
\texttt{TMAC} on the \texttt{fbs} scheme and some user specified parameters (\texttt{params}).
\begin{lstlisting}[caption={example code},label=code:l1log,language=C++]
  // [...] parameters are defined above
  // forward operator: gradient step for logistic loss
  forward_grad_for_log_loss<SpMat> forward(&A, &b, &Atx, eta);
  // backward operator: proximal operator for l1 norm
  prox_l1 backward(eta, lambda);
  // forward-backward splitting scheme
  ForwardBackwardSplitting<forward_grad_for_log_loss<SpMat>, prox_l1>
    fbs(&x, forward, backward);
  // the multicore driver
  TMAC(fbs, params);
\end{lstlisting}
One can easily adapt the previous code to solve other problems, for example, replacing lines 5 through 7 with
the following two lines
\begin{lstlisting}[ caption={modified line 5-7},label=code:tikhonovlog, language=C++]
  prox_sum_square backward(eta, lambda);
  ForwardBackwardSplitting<forward_grad_for_log_loss<SpMat>, prox_sum_square> fbs(&x, forward, backward);
\end{lstlisting}
solves the Tikhonov regularized logistic regression, i.e.,
$$\Min_{x \in \RR^n} \lambda \|x\|_2^2 + \sum_{i = 1}^m \log\left(1 + \exp(-b_i \cdot a_i^T x)\right).$$
Replacing line 3 and line 7 with
the following two lines
\begin{lstlisting}[caption={replaced line 3 and line 7},label=code:LASSO, language=C++]
  forward_grad_for_square_loss<Matrix> forward(&A, &b, &Atx, eta);
  ForwardBackwardSplitting<forward_grad_for_square_loss<Matrix>, prox_l1> fbs(&x, forward, backward);
\end{lstlisting}
solves the Lasso problem
$$\Min_{x \in \RR^n} \lambda \|x\|_1 + \frac{1}{2}\|A x - b\|^2.$$
where \texttt{A} is a dense matrix. It is worth mentioning that the previously used operators have implementions in \pkg. Users can refer to the documentation for
the complete list of implemented operators.

\section{Architecture}

Writing an efficient code is very different from writing an optimization algorithm.
Our toolbox's architecture is designed to mimic how a scientist writes down an optimization algorithm.
The toolbox achieves this by separating into the following layers: Numerical linear algebra, Operator, Scheme, Kernel, and Driver.
Each layer represents a different mathematical component of an optimization algorithm.

The following is a brief description of each layer and how it interacts with the layers above and below it. 

\subsection{Numerical linear algebra}

We use Eigen~\citep{eigenweb}, Sparse BLAS\footnote{\url{http://math.nist.gov/spblas/}}, and BLAS\footnote{\url{http://www.netlib.org/blas/}} in our Toolbox.
Directly using efficient numerical packages like BLAS can be intimidating due to their complex APIs. We provide simplified APIs for common linear algebra operations like $a_i^T x$ in Algorithm \ref{alg:fbs_l1_log}. This layer insulates the user from the grit of raw numerical implementation. Higher layers use the Numerical Linear Algebra Layer in their implementations.
If our provided functions are not sufficient, Eigen, Sparse BLAS, and BLAS are well documented. 

\subsection{Operator}

The Operator Layer contains Forward Operator objects\footnote{``Operator object'' refers to a C++ object. An ``operator'' refers to the mathematical object.} (e.g., gradient descent step, subgradient step) and Backward Operator (e.g., proximal mapping, projection) objects. These operators types see heavy reuse throughtout optimization.
For instance, Algorithm \ref{alg:fbs_l1_log}, along with any other gradient based method for $\ell_1$ regularized logistic regression, requires the computation of the Forward Operator $x - \eta \, \nabla_x \,(\sum_{i = 1}^m \log (1 + e^{-b_i \cdot a_i^T x}))$. On a similar vein, Nonnegative Matrix Factorization and Nonnegative Least Squares share a backwad operator, the projection onto the positive orthant.

Much as the Numerical Linear Algebra Layer insulates the user from the computation details of numerical linear algebra, the Operator Layer insulates the user from the computational details of operators. This is achieved by encapsulating the computation of common forward and backward operators into Operator objects. Abstracting operators into Operator objects is useful, as Operator objects provide clarity and modularity. Consider the construction of  \texttt{forward} in Listing \ref{code:l1log}. At a glance, it is clear that we are using the gradient of logistic regression, the data is sparse, and we are using the Matrix \texttt{A} and Vector \texttt{b} to compute the gradient. Listing \ref{code:l1log} and Listing \ref{code:LASSO} demonstrates how Operator objects provide modularity. Perturbations to an algorithm, such as a change of regularizer or a change of data fidelity term, can be handled by changing the corresponding operator type. The rest of the code structure is unchanged. The Operator Layer, as a result, allows users to reason and code at the level of operators. 
The higher layers use Operator objects as components in the creation of algorithms.
 
As our Toolbox is designed for coordinate update methods, each operator is implemented to compute coordinates efficiently.
 A coordinate or block of coordinates can be computed efficienty if the computational cost of a single coordinate or a block of coordinates of the operator is reduced by a dimensional factor compared to the evaluation of the entire operator (e.g. $(Ax)_i$ versus $Ax$). 
 In some cases caching, the storing an intermediate computations, can improve the efficiency of updating a coordinate block. These ideas are formalized in the paper on coordinate friendly stuctures~\citep{PengWuXuYanYin2016_coordinate}.

In the case that an operator is needed that we have not provided, a user can use~\citep{PengWuXuYanYin2016_coordinate} as a guide to identify coordinate-friendly structures and implement their own operator. In addition, there are certain rules about operator implementation that must be followed; see Section \ref{sc:interface}.

\subsection{Scheme}
A scheme describes how to make a single-iteration update to $x$.
It can be written as a combination of operators. For example, \eqref{eq:fbs_l1_log} is the Forward-Backward Splitting Scheme (also referred as Proximal Gradient Method) for a specific problem, $\ell_1$ regularized logistic regression. In Algorithm \ref{alg:fbs_l1_log}, it corresponds to Line 9. To apply the Forward-Backward Splitting Scheme to the $\ell_1$ regularized logistic regression problem \ref{code:tikhonovlog} , we need to specify a forward operator and a backward operator (e.g., on Line 5 of Algorithm \ref{alg:fbs_l1_log}). We can see this in Listing \ref{code:l1log}.
The scheme object \texttt{fbs} is specialized to $\ell_1$ regularized logistic regression by specifying its type as \texttt{ForwardBackwardSplitting<forward\_grad\_for\_log\_loss<SpMat>, prox\_l1>}.

We provide implemenations of the following schemes:  Proximal-Point Method,  Gradient Descent, Forward-Backward Splitting, Backward-Forward Splitting, Peaceman-Rachford Splitting, and Douglas-Rachford Splitting.

If the provided schemes are not sufficient, the user may implement their own scheme following certain rules so that their scheme can interact with the rest of the package; see Section \ref{sc:interface}.
The user is encouraged to use objects from the Operator Layer as building blocks, but direct calling the Numerical Linear Algebra Layer is perfectly functional. 

\subsection{Kernel}
A Kernel function described the operations that an agent performs.
As can be seen in Algorithm \ref{alg:fbs_l1_log}, it corresponds to the while-loop, which contains contains a coordinate choice rule and a scheme object.
The agent chooses a coordinate using its rule and call the Scheme object with the chosen coordinate to update $x$.
For each coordinate choice rule, we have implemented a corresponding Kernel function.

We provide the following coordinate choice rules: cyclic, random, and parallel Gauss-Seidel. 
Cyclic update rule divides coordinates into approximately equal-sized blocks, and each agent is assigned a block.
Each agent chooses coordinates cyclicly from within the block. 
For random update rule, each agent randomly chooses a block, then chooses coordinates cyclicly from within that block.
For parallel Gauss-Seidel rule, each agent updates all of the coordinates in a Gauss-Seidel fashion. 
If the user desires a different coordinate rule, they may implement their own Kernel function, following the specifications in Section~\ref{sc:implement}.

\subsection{Driver}

Once the schemes and kernels have been specified, agents are responsible for carrying them out. Agents are generated as C++11 threads.
A Driver creates and manages such agents.
For example, if the user chooses the cyclic coordinate Kernel and ten agents, the Multicore Driver will create ten agents using that kernel. 
The Multicore Driver is called with a Params object that contains parameters such as kernel choice, number of iterations, and step size.
An example of this is found in Listing \ref{code:l1log}, where \texttt{fbs} and \texttt{params} are passed to the \texttt{MOTAC} driver.

Optionally, a Multicore Driver can launch a controller agent to control the other agents, for example, by dynamically updating step sizes to accelerate convergence.
The current controller agent monitors convergence by periodically computing the fixed point residual, a measure that would reduce to 0 when the sequence converges to a fixed point.

Most users can treat this layer as a black box. Only when a new way to generate agents is desired, does the user need to modify this layer.

\section{Implementation details}\label{sc:implement}

\subsection{Interaction between layers} \label{sc:interface}

The Operator, Scheme, and Kernel Layers interact heavily with each other. To formalize interaction between each layer, we introduce Layer Interfaces. A Layer Interface describes guaranteed member functions of objects in a Layer so that other Layers can safely use these member functions to interact. The Layer Interfaces allow for specialization of objects within each layer while still maintaining a uniform means of interaction. Consider the Operator Interface:
\begin{lstlisting}[language=C++,label=Operator_Interface]
class OperatorInterface {
public:
  // compute operator at index
  virtual double operator() (Vector* v, int index = 0)=0;
  // compute operator using val at index
  virtual double operator() (double val, int index = 0)=0;
  // compute full operator using v_in, storing in v_out
  virtual void operator() (Vector* v_in, Vector* v_out)=0;
  // update operator related step size
  virtual void update_step_size(double step_size_)=0;
  // update cache variable following an update at index i
  virtual void update_cache_vars(double old_x_i, double new_x_i, int i)=0;
  // update block of cache variables based upon rank of calling thread
  virtual void update_cache_vars(Vector* x, int rank, int num_threads)=0;
};
\end{lstlisting}

For an object to belong to the Operator layer, it must inherit from the Operator Interface.
Code that contains an object that inherits from the Operator Interface that does not implement the functions in the Operator Interface will not compile.

The Scheme Interface:
\begin{lstlisting}[language=C++]
class SchemeInterface {
public:
  // update scheme internal parameters
  virtual void update_params(Params* params)=0;
  // async: compute and apply coordinate update, return S_{index}
  virtual double operator() (int index)=0;
  // sync: compute and store S_{index} in S_i
  virtual void operator() (int index, double& S_i)=0;
  // sync: apply block of S stored in s to solution variable
  virtual void update(Vector& s, int range_start, int num_cords)=0;
  // sync: apply coordinate of S stored in s to solution variable
  virtual void update(double s, int idx)=0;
  // sync: update rank worth of cache_vars based on num_threads
  virtual void update_cache_vars(int rank, int num_threads)=0;
};
\end{lstlisting}
For an object to belong to the Scheme Layer, it must inherit from the Scheme Interface.
Code that contains an object that inherits from the Scheme Interface that does not implement the functions in the Scheme Interface will not compile.

Schemes must be suitable for both asynchronous and synchronous computing. Synchronous computing requires all coordinate updates to be computed before applying the coordinate updates. Asynchronous computing can apply coordinates updates immediately. The Scheme Interface reflects these two computing regimes.

The Scheme Interface is very lightweight, as the iterations of optimization methods come in a variety of forms.  All it requires is that a coordinate update can be produced and that parameters can be passed to the object.

\subsection{Kernel and Multicore Driver interaction}

A Multicore Driver creates agents using Kernel functions. The arguments to a Multicore Driver are a scheme object and a params object. A new Kernel function must be callable using information from those two objects. The Multicore Drivers must be modified to include the new kernel function as an option for creating an agent. The modfication requires adding a new case to an if-elseif chain, and knowledge of how to create a C++11 thread.

\subsection{Templating}

In C++, when objects have similar structures that only vary based upon an input type, templates are used to reduce code redundancy. For instance, the code for an object representing a dense matrix of \emph{doubles} and the code for an object representing a dense matrix of \emph{floats} is identical.
Templates are not objects, but instead are blueprints for constructing an object.
Based upon the arguments to the template, a corresponding type is automically constructed.
We use templating heavily in our Toolbox, as most of our workflow can be genericized.

Objects in the Scheme Layer are implemented as templates.
This is a natural choice, as the Scheme objects are generic iteration formulas.
For example, in Listing \ref{code:l1log}, the scheme type of \texttt{fbs} is defined by the arguments to the forward backward splitting template, \texttt{forward\_grad\_for\_log\_loss<SpMat>}, and \texttt{prox\_l1}.

Objects in the Object Layer are also templatized. Depending on data representation, it can be more efficient to use functionality designed for that data representation. Consider the difference between computing $x - \eta \, \nabla_x \,(\sum_{i = 1}^m \log (1 + e^{-b_i \cdot a_i^T x}))$ when $a_i$ is stored in the  sparse format\footnote{Only the values and locations of  nonzero entries are stored.} instead of the dense format.
Our linear algebra functions are overloaded so the compiler will deduce the proper function to use in the template.

Kernel functions are also templatized. This allows Kernel functions to take in as input arbitrary objects from the Scheme Layer. If we did not use templating, for each coordinate rule would need a function for every possible realization of gradient descent, proximal point method, etc.  Any object that satisfies the Scheme Interface can be passed to a Kernel function.

Multicore Drivers are templatized for the same reason as Kernel functions are templated. If we did not use templating, we would need a Multicore Driver for every possible realization of gradient descent, proximal point method, etc. Any object that satisfies the Scheme Interface can be passed to a Mutlicore Driver.






\section{Numerical experiments}
In this section, we illustrate the efficiency of \pkg~for three applications: $\ell_1$ regularized logistic regression, portfolio optimization,  and nonnegative matrix factorization. The tests were carried out on a machine with 64GB of memory and two Intel Xeon E5-2690 v2 processors (20 cores in total).

\subsection{Minimizing $\ell_1$ logistic regression}
In this subsection, we apply \pkg~to the $\ell_1$ regularized logistic regression problem~\eqref{eq:l1_log}. It implements Algorithm~\ref{alg:fbs_l1_log}. The command-line executable is tmac\_fbs\_l1\_log. We set $\lambda =0.001$, and tests TMAC on two LIBSVM datasets\footnote{\url{http://www.csie.ntu.edu.tw/~cjlin/libsvmtools/datasets/}}: news20, and url.

\begin{figure}[!h]
        \centering
       \begin{subfigure}[b]{0.4\textwidth}
                \includegraphics[width=\textwidth]{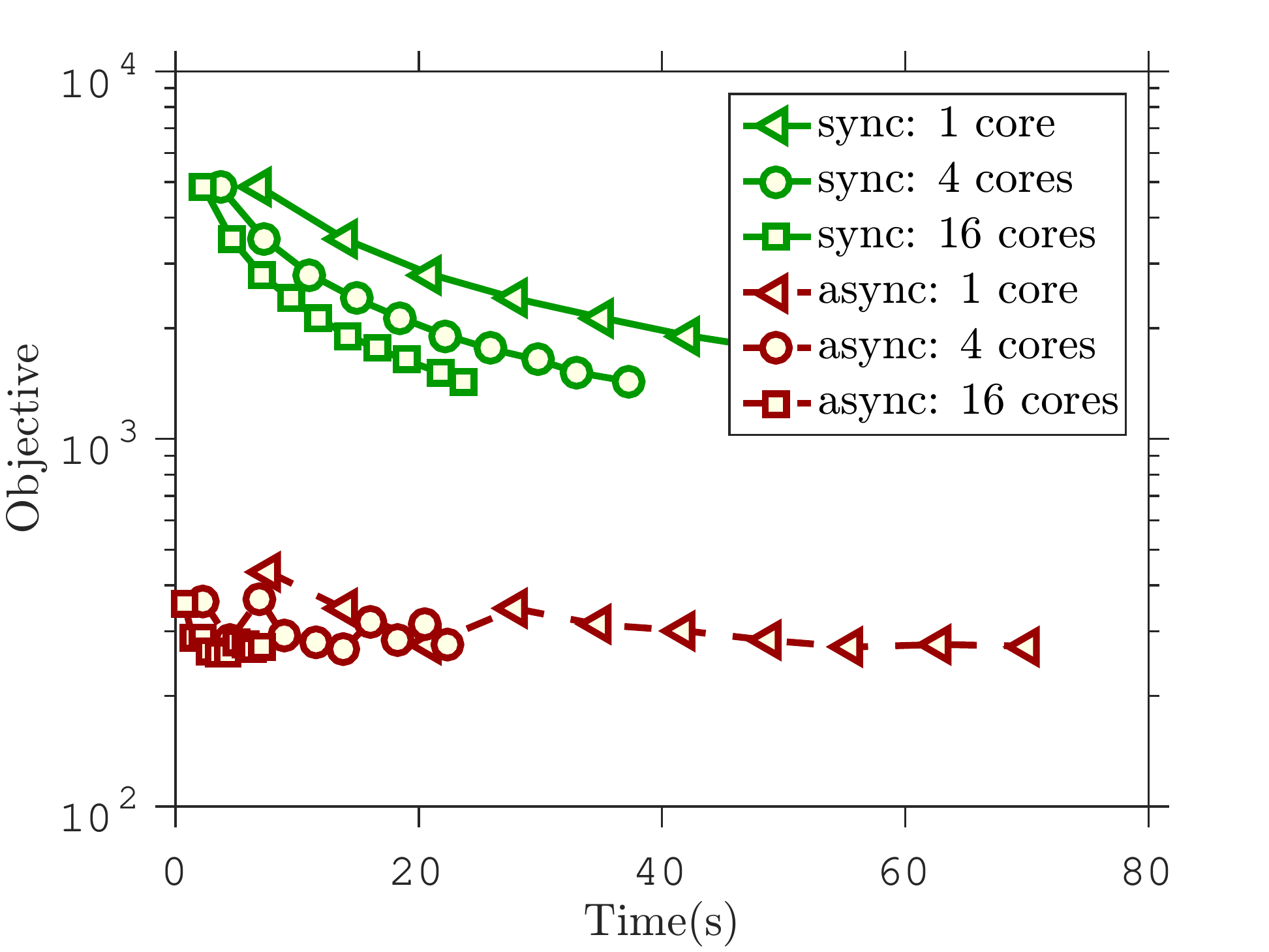}
                \caption{dataset ``news20''}
        \end{subfigure}
        ~~
        \begin{subfigure}[b]{0.4\textwidth}
                \includegraphics[width=\textwidth]{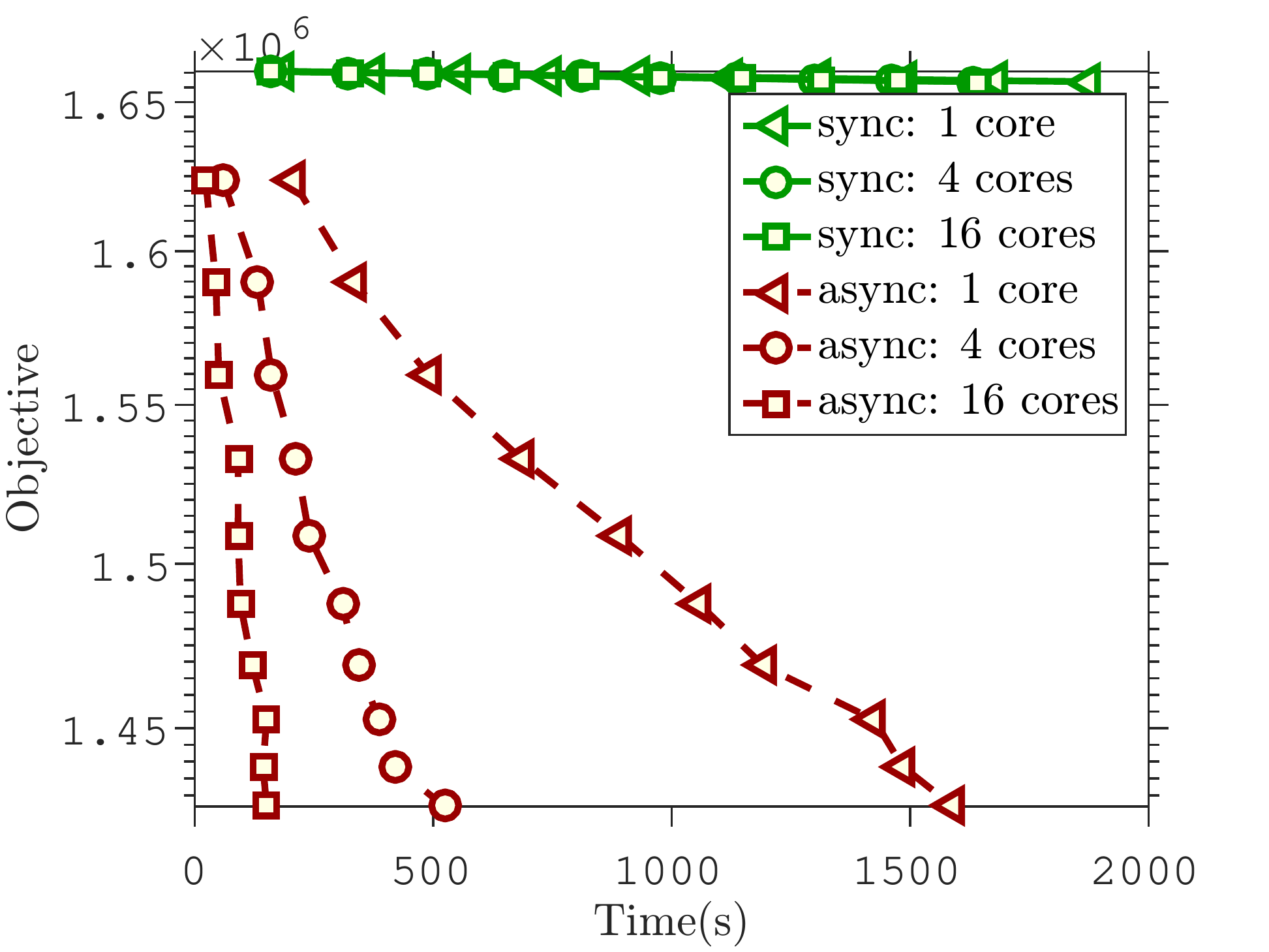}
                \caption{dataset ``url''}
        \end{subfigure} 
        \caption{Objective vs wall clock time.}\label{fig:log_reg_obj}
\end{figure}

\begin{figure}[!h]
        \centering
       \begin{subfigure}[b]{0.35\textwidth}
                \includegraphics[width=\textwidth]{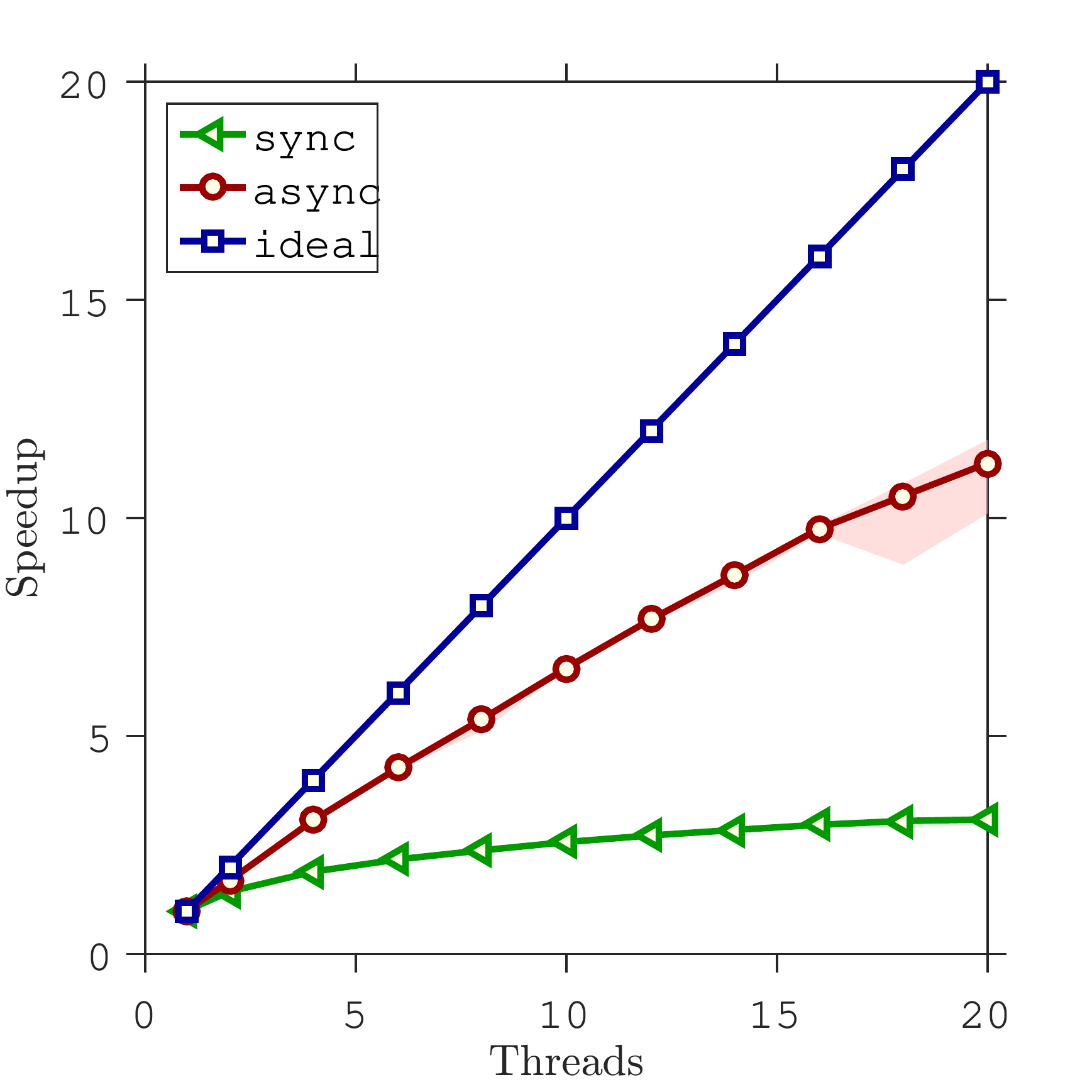}
                \caption{dataset ``news20''}
        \end{subfigure}
        ~~
        \begin{subfigure}[b]{0.35\textwidth}
                \includegraphics[width=\textwidth]{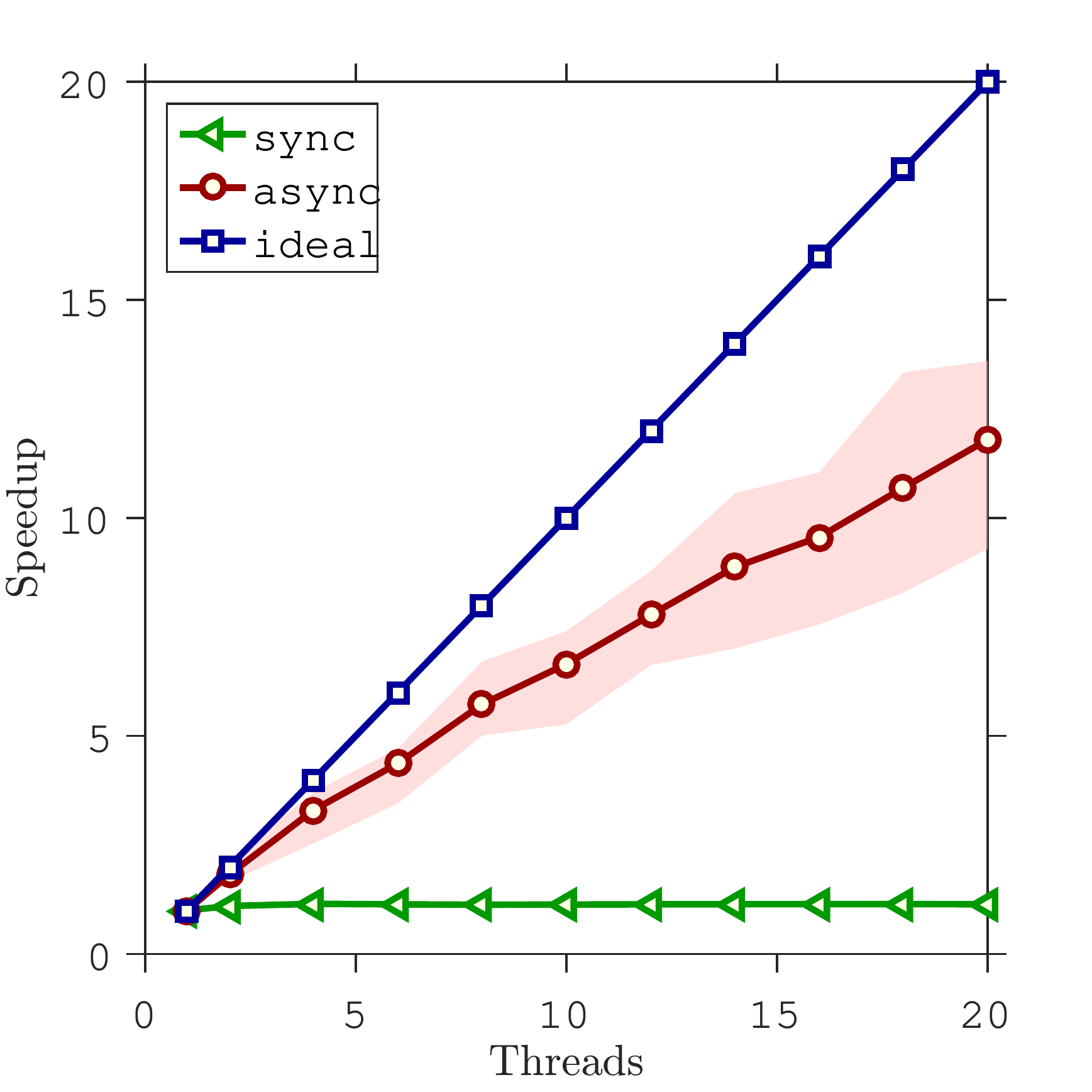}
                \caption{dataset ``url''}
        \end{subfigure}        
        \caption{Speedup vs number of threads. The solid lines represent mean speedup across 10 different runs. The shaded regions represent the lower and upper bounds of speedup for the 10 runs.}\label{fig:log_reg_speedup}
\end{figure}

Figure \ref{fig:log_reg_obj} gives the running times of  the sync-parallel and async-parallel implementations on the two datasets. Figure~\ref{fig:log_reg_speedup} is the speedup performance comparison of the two methods. We can observe that async-parallel achieves near-linear speedup, but sync-parallel scales poorly as we shall explain below. One can also see that async-parallel converges faster due to more relaxed forward operator step size selection.

In the sync-parallel implementation,  all the running cores have to wait for the last core to finish an iteration, and therefore if a core has a large load, it slows down the iteration. Although every core is (randomly) assigned to roughly the same number of features at each iteration, their  $a_i$'s have very different numbers of nonzeros, and the core with the largest number of nonzeros is the slowest (Sparse matrix computation is used for both datasets, which are very large.) As more cores are used,  despite that they altogether do more work at each iteration, the per-iteration time reduces as the slowest core tends to be relatively slower.


\subsection{Portfolio optimization}
Assume that we have one unit of capital and $m$ assets to invest on. The $i$th asset has an expected return rate $\xi_i\ge 0$\cut{ for $i\in[m]$}. Our goal is to find a portfolio with the minimal risk such that the expected return is no less than $c$. This problem can be formulated as

\begin{equation*}
\begin{array}{l}
\displaystyle
\Min_x ~ \frac{1}{2} x^\top Q x, \\
\displaystyle
\text{subject to}~ x\ge0, \sum_{i=1}^m x_i\le 1,\, \sum_{i=1}^m\xi_i x_i\ge c,
\end{array}
\end{equation*}
where the objective function is a measure of risk, and the last constraint imposes that the expected return is at least $c$. We apply the coordinate update scheme as shown in Section 5.3.1 of~\citep{PengWuXuYanYin2016_coordinate} to solve it.
We tested two synthetic  problem instances: one has 5,000 assets and the other instance has 30,000 assets.
The vector of expected return rate was sampled from $N(0.01, 1)$ normal distribution.
The matrix $Q$ was set to $\frac{1}{2} (R + R^T) + \sigma \cdot I$, where the entries of $R$ was sampled independently from $N(0, 0.1)$ normal distribution, and $\sigma$ was chosen such that $Q$ was positive definite.
We tested both sync-parallel and async-parallel methods with 100 epochs. They reached similar final objective. We report the timing and speedup results in Table~\ref{tab:port_opt}. One can observe that  TMAC scales well for both sync-parallel method and async-parallel method. The nice scaling performance of sync-parallel method is due to almost perfect load balancing across the threads and the homogeneous computing environment.
\begin{table}[!h]
\centering
\begin{tabular}{rrrrrr}
\toprule
\multirow{2}{*}{} & &  \multicolumn{2}{c}{async-parallel} & \multicolumn{2}{c}{sync-parallel} \\
\midrule
 & & time (s) & speedup & time (s) & speedup \\
 \midrule
 \multirow{3}{*}{5K assets} & 1 thread~ & 13.95& 1.00 & 13.98& 1.00 \\
  & 4 threads & 3.53& 3.95 &3.54 & 3.95\\
 & 16 threads & 0.99 &14.09 &0.92 & 15.19\\
\midrule
 \multirow{3}{*}{30K assets} & 1 thread~ & 483.20& 1.00& 479.22& 1.00 \\
  & 4 threads & 124.94& 3.86 &123.74 & 3.87\\
 & 16 threads & 31.96 &15.12 &33.11 & 14.47\\
\bottomrule
\end{tabular}
 \caption{\label{tab:port_opt}Results for portfolio optimization. }
\end{table}

\subsection{Nonnegative matrix factorization}

We consider the following Nonnegative Matrix Factorization problem

\begin{equation*}
	\Min_{X \geq 0,Y \geq 0} ~\|A-X^T Y\|_2^2,
\end{equation*}
where $A \in \RR^{m\times n}$, $X \in \RR^{k \times m}$ and $Y \in \RR^{k \times m}$.
This problem, despite being nonconvex, has a special form.
The objective function is block multiconvex\footnote{the objective function is convex when all but a few specific variables are held fixed} and its reguralizers are separable. Problems of this type have been shown \citep{XuYin2013_block,BolteSabachTeboulle2014_proximal} to be amenable to coordinate update techniques.
Recent work\citep{2016APALM} has shown problem of this type to amenable to the asynchronous regime.
As the problem is nonconvex, convergence is given to a local minimizer, not a global minimizer.

We run TMAC on a synthetic problem, $A=\hat X^T \hat Y$ ,  with $m=1000$ and $k=20$.
Elements of $\hat X$ and $\hat Y$ sampled independently from $N(0, 1)$ normal distribution, then thresholded positive.
We ran the tests with variable number of threads and iterations. 
The following results are the averages resulting from 20 runs.

\begin{figure}[!h]
        \centering
       \begin{subfigure}[b]{0.4\textwidth}
                \includegraphics[width=\textwidth]{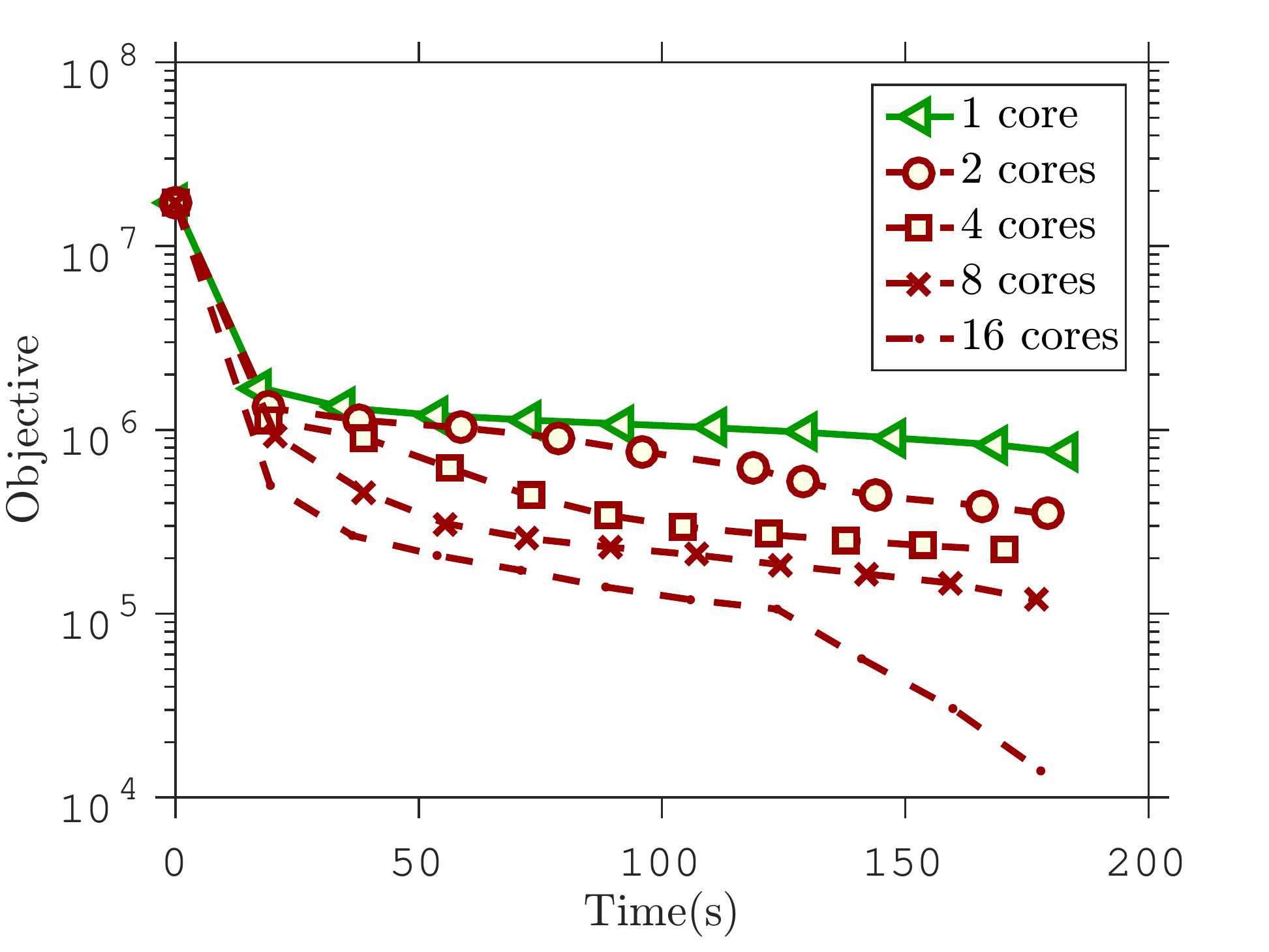}
                \caption{Objective vs walll clock time}
        \end{subfigure}
        ~~
        \begin{subfigure}[b]{0.4\textwidth}
                \includegraphics[width=\textwidth]{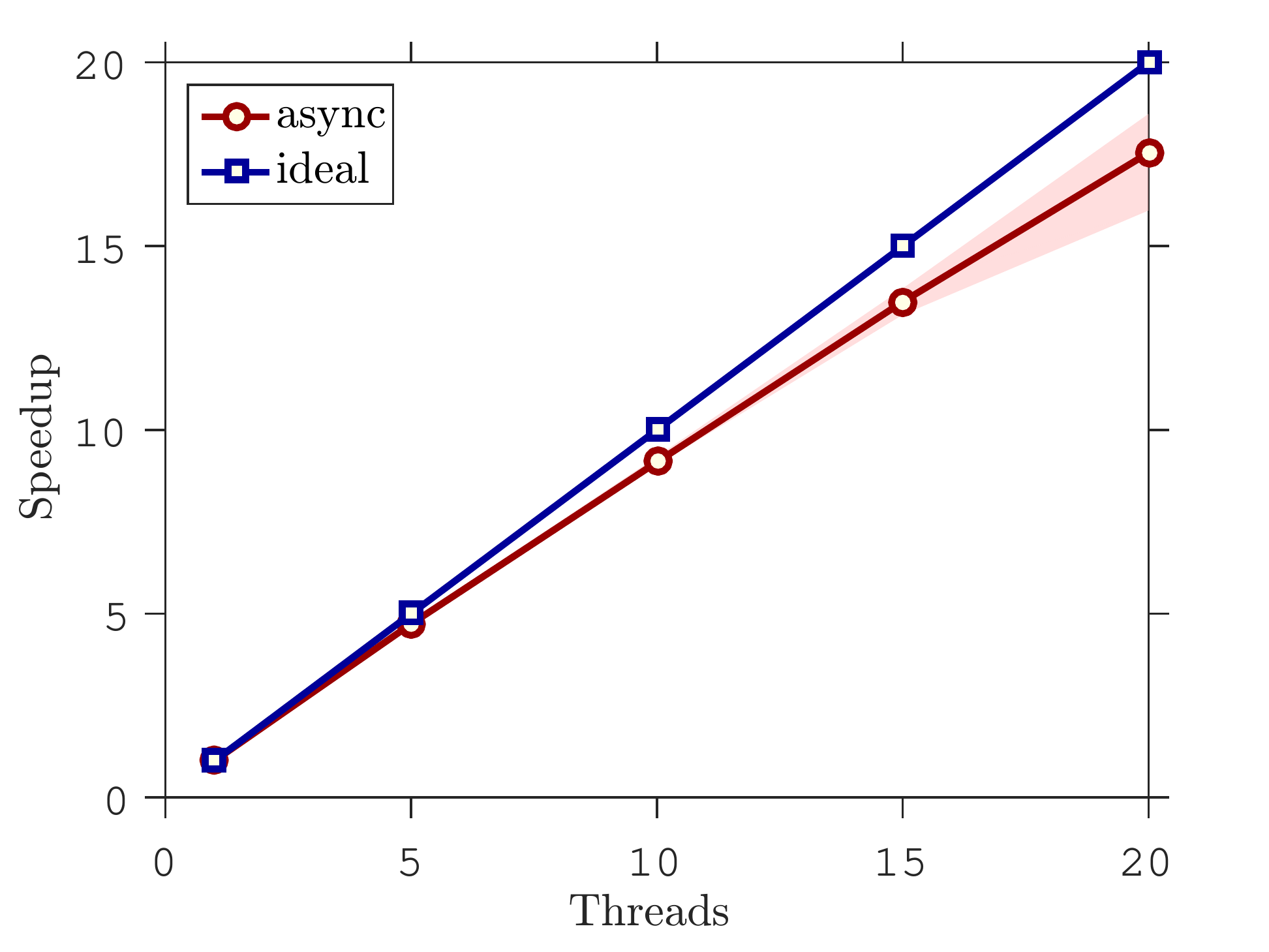}
                \caption{speedup}
        \end{subfigure} 
        \caption{NMF results}\label{fig:NMF_speedup}
\end{figure}

To show scalability we increased the dimension of $k$, and tested the speedup.
\begin{table}[!h]
\centering
\begin{tabular}{rrrr}
 \toprule
{threads}&{k=10}&{k=20}&{k=100}\\\cmidrule{1-4}
1&1.00&1.00&1.00\\
2&1.97&1.98&1.98\\
4&3.75&3.75&3.76\\
8&7.12&7.33&7.35\\
16&13.38&14.51&14.43\\
\bottomrule
\end{tabular}
 \caption{\label{tab:nmf_res}Speedup results for nonnegative matrix factorization. }                  
\end{table}



\section{Future Work}

New features are still being added to \pkg. The following are some of the features we are exploring.

\subsection{Stochastic algorithms} 

Stochastic algorithms exploit summative structures in problems to produce cheap updates. Our current toolbox does not currently support stochastic algorithms, but can be modified to do so. Such a modification would require stochastic operators. A developement branch will appear on our github.  

\subsection{Cluster computing}
 
 Currently, agents are realized as threads.
 This limits our toolbox to the multicore regime.
 Future releases intend to use MPI to bring \pkg~to cluster computing. The new functionality will be provided in the Driver and the Kernel Layer. 

\subsection{User interface}
 \pkg~requires the user to either use our prepackaged executables, or code moderately in C++.
 A graphic user interace is in development for those who wish avoid interacting directly with C++.
 This will limit the user to built-in functionality.
 In addition, interfaces to Matlab and Python will be provided for algorithms of particular interest.
 
\subsection{Automatic parameter choice}
\pkg~requires the user to choose stepsizes, and number of threads.
In future releases we intend to provide automatic stepsize heuristics. 
Optimizing thread number is a function of current processor usage and computing architecture. 
We intend to provide functions to survey the current architecture and suggest proper levels of parallelism and asynchrony.

\subsection{Block coordinate update}
  Currently \pkg~does not support block coordinate updates (updates consisting of a set of coordinates). Block coordinate updates present a tradeoff between iteration complexity and communication efficiency.
  Automatic block size deduction and block composition (the coordinates forming a block) is an open question.
We plan to explore several heuristics.

\subsection{New fields}

Splitting methods have been used in many fields (see \citep{roland2016some} for a more in depth discussion).
Our current release focuses on optimization, but provides a strong foundation for branching into other fields.
Fruitful fields to explore include:

\begin{itemize}
\item Numerical simulations;
\item Large scale numerical linear algebra;
\item Time varying systems such as initial value problems and partial differential equations.
\end{itemize}

\section{Conclusion }
We have developed \pkg, an easy-to-use open source toolbox for large scale optimization problems.
The toolbox implemented both sequential and parallel algorithms based on operator splitting methods, stochastic methods,
and coordinate update methods. \pkg~is separated into several layers and mimics how a scientist writes down an optimization algorithm. Therefore, it is easy for one to obtain a new algorithm by modifying just one of the layers such as adding a new operator.

New features and applications will be added to \pkg~based upon new research and community input. The software and user guide \repo~will be kept up-to-date and supported.

\section{Acknowledgements}
The work is supported in part by NSF grants ECCS-1462398.

\bibliography{wotaoyin_lib,asyn}

\end{document}